\documentclass{icmart}
\usepackage{graphicx}

\contact[francisco.santos@unican.es, \ URL: http://personales.unican.es/santosf/]{Francisco Santos,
Departamento de Matem\'aticas, Estad\'{\i}stica y Computaci\'on, Universidad de Cantabria, E-39005 Santander, SPAIN}

\newcommand{\invlim}{\varprojlim}
\newcommand{\link}{{\operatorname{link}}}

\newcommand{\reals}{\mathbb{R}}
\newcommand{\integers}{\mathbb{Z}}
\newcommand{\naturals}{\mathbb{N}}
\newcommand{\T}{{T}}
\newcommand{\A}{\mathcal{A}}
\newcommand{\B}{\mathcal{B}}
\newcommand{\M}{\mathcal{M}}
\newcommand{\vertices}{\operatorname{vertices}}
\newcommand{\conv}{\operatorname{conv}}
\newcommand{\pos}{\operatorname{pos}}

\newtheorem{theorem}{Theorem}[section]
\newtheorem{corollary}[theorem]{Corollary}
\newtheorem{lemma}[theorem]{Lemma}
\newtheorem{proposition}[theorem]{Proposition}
\newtheorem{conjecture}[theorem]{Conjecture}

\newtheorem*{thm}{Theorem}
\newtheorem*{cor}{Corollary}

\theoremstyle{definition}
\newtheorem{definition}[theorem]{Definition}
\newtheorem{remark}[theorem]{Remark}
\newtheorem{ejemplo}[theorem]{\small Example}
\newenvironment{example}
{\begin{ejemplo}\footnotesize
}
{\end{ejemplo}}

\title[Geometric bistellar flips]
{Geometric bistellar flips.
The setting, the context and a construction}

\author[Francisco Santos]{Francisco Santos
\thanks{Partially supported by  the Spanish Ministry of Education and Science, grant number
 MTM2005-08618-C02-02.}
}

\begin{document}

\begin{abstract}
We give a self-contained introduction to the theory
of secondary polytopes and geometric bistellar flips
in triangulations of polytopes and point sets,
as well as a review of some of the known results and
connections to algebraic geometry, topological combinatorics,
and other areas.

As a new result, we announce the construction of a point set
in general position with 
a disconnected space of triangulations. This shows, for the first time, that the
poset of strict polyhedral subdivisions of a point set
is not always connected.
\end{abstract}



\maketitle

\section*{Introduction}
\label{sec:intro}

Geometric bistellar flips are ``elementary moves'', that is,
minimal changes, between triangulations of a point set
in affine space $\reals^d$. In their present form
they were introduced around 1990 by
Gel'fand, Kapranov and Zelevinskii during their study
of discriminants and resultants for sparse polynomials~\cite{GKZ-paper,GKZ-book}.
Not surprisingly, then, these bistellar flips have several
connections to algebraic geometry. For example,
the author's previous constructions of point sets 
with a disconnected graph of 
triangulations in dimensions five and 
six~\cite{Santos-noflips,Santos-hilbert}
imply that certain algebraic schemes considered
in the literature~\cite{Alexeev,BiKaSt,HaiStu,PeeSti}, including the so-called
toric Hilbert scheme, are sometimes not connected.

Triangulations of point sets play also an obvious role in
applied areas such as computational geometry or computer
aided geometric design,
where a region of the plane or 3-space is triangulated in order to
approximate a surface, answer proximity or visibility questions, 
etc. See, for example, the survey articles~\cite{AurenKlein,Bern-survey}, or ~\cite{Edels-book}. In these fields, flips
between triangulations have also been considered since 
long~\cite{Lawson}. Among other things,
they are used as the basic step to compute an optimal
triangulation  of
a point set incrementally, that is, adding the points one by one.
This \emph{incremental flipping algorithm} is the one usually preferred for,
for example, computing the Delaunay triangulation,
as ``the most intuitive and easy to implement''~\cite{AurenKlein},
and yet as efficient as any other.

%

In both the applied and the theoretical framework, 
the situation is the same: a fixed set of points
$\A\subset\reals^d$ 
is given to us (the ``sites'' for a Delaunay triangulation computation,
the test points for a surface reconstruction, or a set of monomials,
represented as points in $\integers^d$, in the
algebro-geometric context) and we need to either explore the collection of
all possible triangulations of this set 
$\A$ or search for a particular one
that satisfies certain optimality properties. Geometric bistellar
flips are the natural way to do this. For this reason, it
was considered one of the main open questions in polytope theory 
ten years ago whether point sets exist with triangulations that cannot be connected via these 
flips~\cite{Ziegler-recent}. As we have mentioned above, 
this question was answered negatively by the author of this
paper, starting in dimension 5. The question is still open in dimensions three and four.

\medskip
This paper intends to be an introduction to this topic, organized in three parts. 

The first section is a self-contained introduction to the theory of
geometric bistellar flips and \emph{secondary polytopes} in triangulations
of point sets, aimed at the non-expert. The results in it are
certainly not new (most come from the original work of Gel'fand, Kapranov
and Zelevinskii mentioned above)
but the author wants to think that this section has some expository novelty; several examples that illustrate the theory 
are given, and our introduction of geometric bistellar flips 
first as
certain polyhedral subdivisions and only afterwards as transformations
between triangulations is designed to show that the definition
is as natural as can be. 
This section 
finishes with an account of the state-of-the-art
regarding knowledge of the graph of flips for sets with ``few'' points
or ``small'' dimension, with an
emphasis on the differences between dimensions two and three.

The second section develops in more detail the two contexts in which we have mentioned that flips are interesting
(\emph{computational geometry} and \emph{algebraic geometry}) together with 
other two, that we call ``\emph{combinatorial topology}'' and
``\emph{topological combinatorics}''. Combinatorial topology refers to the
study of topological manifolds via triangulations of them.
Bistellar flips have been proposed
as a tool for manifold recognition~\cite{BjornerLutz,Lickorish},
and triangulations of the 3-sphere without bistellar flips other
than ``insertion of new vertices''
are known~\cite{DoFaMu}. Topological combinatorics 
refers to topological methods
in combinatorics, particularly to the topology of 
\emph{partially
ordered sets} (posets) via their order complexes. The graph of 
triangulations of a point set $\A$ consists of the first two levels
in the poset of polyhedral subdivisions of $\A$, which in turn
is just an instance of several similar posets studied in combinatorics
with motivations and applications ranging from oriented matroid theory
to bundle theories in differential geometry.

The third  section announces for the first time the construction
of a point set \emph{in general position} whose graph of triangulations
is not connected. The details of the proof appear in~\cite{Santos-17}. The point set is also the smallest one known so far
to have a disconnected graph of flips.

\begin{thm}
\label{thm:intro}
There is a set of 17 points in general position in $\reals^6$ whose
graph of triangulations is not connected.
\end{thm}

As usual in geometric combinatorics, 
a finite point set $\A\subset\reals^d$ is said to be \emph{in general 
position} if no $d+2$ of the points lie in an affine hyperplane. Equivalently,
if none of the $\genfrac{(}{)}{0pt}{1}{|\A|}{d+1}$ determinants defined by the point set
vanish. Point sets in general position form an open dense subset in the space $\reals^{n\times d}$ of sets of dimension $d$ with $n$ elements. That is to say,
``random point sets'' are in general position. Point sets that are not in general position are said to be \emph{in special position}.
\medskip

The connectivity question  
has received special attention  in general position
even before disconnected examples in special position were found.
For example, Challenge 3 in \cite{Ziegler-recent} and Problem 28 in \cite{CG-open}  specifically ask
whether disconnected graphs of flips exist for point sets in special position
(the latter asks this only for dimension 3). 
Although it was clear (at least
to the author of this paper) from the previous examples 
of disconnected graphs of flips  that
examples in general position should also exist, modifying those particular
examples to general position and proving that their flip-graphs are
still not connected is
not an easy task for quite intrinsic reasons: the proofs of non-connectednes
in~\cite{Santos-noflips,Santos-hilbert} 
are based on the fact that the point sets considered 
there are cartesian products of 
lower dimensional ones.

In our opinion, 
an example of a disconnected graph of flips
in general position is interesting  for the following three reasons:

\begin{enumerate}

\item The definition of flip that is most common in computational geometry coincides with  ours (which is the standard one in algebraic geometry
and polytope combinatorics) only for point sets in general position. In special position, the computational geometric definition is far more restrictive and, in particular, taking it
makes disconnected graphs of flips in special
position be ``no surprise''. For example, Edelsbrunner~\cite{Edels-book}
says that the flip-graph among the (three) triangulations of a regular octahedron
is not connected; see Section~\ref{sub:computational-geometry}.

\item Leaving aside the question of definition, in
engineering applications the coordinates of points are usually 
approximate and there is no loss in perturbing them  
into general position. That is, the general position case is sometimes
the only case.

\item Even in a purely theoretical framework, point sets in general position
have somehow simpler properties
than those in special position. If a point set $\A$
in special position has a non-connected graph of flips then
 automatically some subset of $\A$
(perhaps $\A$ itself) has a disconnected poset of subdivisions.
This poset
is sometimes called the \emph{Baues poset} of $\A$ and its study is (part of)
the so-called generalized Baues problem. See Section~\ref{sub:topological-combinatorics}, or~\cite{Reiner-survey} for more
precise information on this. 
In partiular, the present example is the first one
(proven) to have a disconnected Baues poset.
\end{enumerate}

\begin{cor}
There is a set of at most 17 points in $\reals^6$ whose poset of proper
polyhedral subdivisions is not connected.
\end{cor}


\section{The setting}
\label{sec:setting}

\subsection{Triangulations. Regular triangulations and subdivisions}
\label{sub:regular}

\subsubsection*{Triangulations and polyhedral subdivisions}
A (convex) \emph{polytope} $P$ is the convex hull of a finite set of points in 
the affine space $\reals^d$.
A \emph{face} of  $P$ is its intersection with any hyperplane 
that does not cross the relative interior of $P$. (Here, the 
\emph{relative interior}
 of  $S\subseteq\reals^d$ is the interior of $S$ regarded as 
a subset of its affine span).
We remind the reader that the faces of dimensions $0,1,d-2$ and $d-1$ of
a $d$-polytope are called vertices, edges, ridges and facets, respectively.
Vertices of $P$ form the minimal $S$ such that $P=\conv(S)$.

A $k$-simplex is a polytope whose vertices (necessarily $k+1$)
are affinely independent. 
It has $\genfrac{(}{)}{0pt}{1}{k+1 }{ i+1}$
faces of each dimension $i=0,\dots,k$, which are all simplices.

\medskip

\begin{definition}
Let $\A$ be a finite point set in $\reals^d$.
\label{defi:triangulation}
A {\em triangulation} of $\A$ is any  collection $\T$ of 
affinely spanning and 
affinely independent
subsets of $\A$ with the following  properties: 
\begin{enumerate}
\item if $\sigma$ and $\sigma'$ are in $\T$, then
$\conv(\sigma)\cap\conv(\sigma')$ is a face of both 
$\conv(\sigma)$ and $\conv(\sigma')$. That is,
$\T$ induces a geometric simplicial complex in $\reals^k$; 
\item $\cup_{\sigma\in
\T}\conv(\sigma)=\conv(\A)$. That is,  $\T$ covers the convex hull of $\A$.
\end{enumerate}
\end{definition}

Note that our definition allows
for some points of $\A$ not to be used at all in a particular triangulation.
Extremal points (vertices of $\conv(\A)$) are used in every triangulation.
%
The elements of a triangulation $T$ are called \emph{cells}.


We can define 
\emph{polyhedral subdivisions} of $\A$ by removing the requirement of the
sets $\sigma$ to be affinely independent in 
Definition~\ref{defi:triangulation}. Since a general subset
$\sigma$ of $\A$ may contain points which are not vertices of $\conv(\sigma)$,
now the fact that the elements of a subdivision are subsets of $\A$
rather than ``subpolytopes'' is not just a formality: 
points which are not vertices of any ``cell''
in the subdivision may still be considered ``used''
as elements of some 
cells. In order to get a nicer concept of polyhedral subdivision,
we also  modify part 1 in 
Definition~\ref{defi:triangulation},
adding the following (redundant for affinely independent sets) condition:
\[
\forall \sigma,\sigma'\in T,\qquad
\conv(\sigma\cap\sigma')\cap \sigma=
\conv(\sigma\cap\sigma')\cap \sigma'.
\]
That is, if $\A$ contains some point
 in the common face $\conv(\sigma\cap\sigma')$ 
of $\conv(\sigma)$ and $\conv(\sigma')$ but not a vertex of it,
that point is either in both or in none of $\sigma$ and $\sigma'$.

Polyhedral subdivisions of $\A$ form a \emph{partially ordered set}
(or \emph{poset})
with respect to the following refinement relation:
\[
S\text{ refines } S' \qquad
:\Leftrightarrow \qquad
\forall \sigma'\in S', \exists \sigma\in S,\text{ such~that } \sigma\subseteq \sigma'.
\]
Triangulations are, of course, the minimal elements in this poset.
The poset has a unique maximal element, 
namely the \emph{trivial suvdivision} $\{\A\}$.

\begin{example}
\label{exm:fivepoints}
Let $\A$ be the following set of five points $a_1,\dots,a_5$
in the plane. We take
the convention that points are displayed as columns in a matrix, 
and that an extra homogenization coordinate (the row of 1's in the following
matrix) is added so that linear algebra, rather than affine geometry,
can be used for computations:
\begin{align}
\label{matrix:fivepoints}
\A=
\bordermatrix{
& a_1 & {a_2} & {a_3} & {a_4} & {a_5}  \cr
&0&3&0&3&1\cr
&0&0&3&3&1\cr
&1&1&1&1&1\cr
}
\end{align}
The following are the nine polyhedral subdivisions of $\A$. Arrows represent the refinement relation, pointing
from the coarser to the finer subdivision. 
For clarity, we write ``$125$'' meaning
$\{a_1,a_2,a_5\}$, and so on.
Figure~\ref{fig:fivepoints} shows pictures of the subdivisions.
In the corners are the four triangulations of $\A$ and  in the middle is the trivial subdivision.
\[
\begin{matrix}
\{125,135,235,234\} & \leftarrow & \{1235,234\} & \to  & \{135,234\}\cr
\uparrow && \uparrow & & \uparrow\cr
\{125,135,2345\} & \leftarrow &\{12345\} & \to & \{1234\}\cr
\downarrow && \downarrow && \ \downarrow\cr
\{125,135,245,345\} & \leftarrow & \{1245,1345\}& \to & \{124,134\}\cr
\end{matrix}
\]
\begin{figure}[htb]
\centerline{
\includegraphics[height=2.5 in]{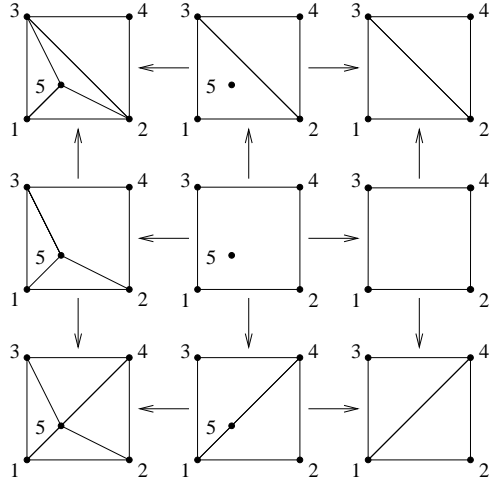}{}{}
}
\caption{The nine polyhedral subdivisions of a certain point set.}
\label{fig:fivepoints}
\end{figure}
The last two columns of subdivisions geometrically induce
the same decomposition of $\conv(\A)$ into subpolygons. Still, we consider
them different subdivisions since the middle column ``uses'' the interior
point $5$ while the right column does not.
\end{example}

\subsubsection*{Regular subdivisions}
Let a point set $\A$ be given, and choose 
a function $w:\A\to \reals$ to lift $\A$ to $\reals^{d+1}$ as the
point set
\[
 \A_w:=\{(a,w(a) : a\in \A\}.
\]
A \emph{lower facet} of $\conv( \A_w)$ is a facet whose supporting hyperplane lies below the interior of  $\conv( \A_w)$. 
The following is a polyhedral subdivision of $\A$, where $\pi:\reals^{d+1}\to \reals^d$ is the projection that forgets the last coordinate:
\[
T_w:=\{\pi(F\cap  \A_w): F \text{ is a lower facet of }\conv( \A_w)\}.
\]
Geometrically, we are projecting down onto $\A$ the lower envelope
of $\A_w$, keeping track of points that lie in the lower boundary 
even if they are not vertices of a facet.

\begin{definition}
\label{defi:regular}
The polyhedral subdivisions and triangulations 
that can be obtained in this way are called
\emph{regular}.
\end{definition}

If $w$ is sufficiently generic then $T_w$ is clearly a triangulation.
Regular triangulations
are  particularly simple and yet quite versatile. 
They appear in different
contexts under different names such as \emph{coherent}~\cite{GKZ-book}, \emph{convex}~\cite{ItenbergShustin,Viro}, \emph{Gale}~\cite{McMullen},
or \emph{generalized (or, weighted)  Delaunay}~\cite{Edels-book}
triangulations.  The latter 
refers to the fact that the Delaunay triangulation of $\A$, probably the
most used triangulation in applications, is the regular 
triangulation obtained with  $w(a)=||a||^2$, where $||\cdot||$
is the euclidean norm.

\begin{example}
\label{exm:moae}
Let
\[
\A=
\bordermatrix{
& {a_1} & {a_2} & {a_3} & {a_4} & {a_5} & {a_6} \cr
&4 & 0 & 0 & 2 & 1 & 1 \cr
&0 & 4 & 0 & 1 & 2 & 1 \cr
&0 & 0 & 4 & 1 & 1 & 2 \cr
}.
\]
This is a configuration of six points in the affine plane with equation
$x_1+x_2+x_3=4$ in $\reals^3$. Since the matrix is already homogeneous
(meaning precisely that columns lie in an affine hyperplane) we do not
need the extra homogenization row. 
The configuration consists of two 
parallel equilateral triangles, one inside the other.
We leave it to the reader to check that the following are two non-regular
triangulations (see Figure~\ref{fig:moae-triangs}):
\[
T_1:=\{124, 235, 136, 245, 356, 146, 456\},\quad
T_2:=\{125, 236, 134, 145, 256, 346, 456\}.
\]
\begin{figure}[htb]
\centerline{
\includegraphics[height=1 in]{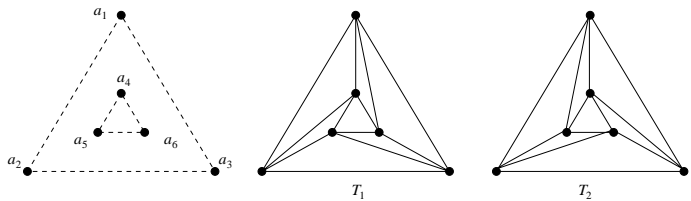}{}{}
}
\caption{A point configuration with two non-regular triangulations}
\label{fig:moae-triangs}
\end{figure}
This example is the smallest possible, since 1-dimensional point configurations
and point configurations with at most $d+3$ points in any dimension $d$ only
have regular triangulations. The former is easy to prove and the latter was first shown in~\cite{Lee-regular}.
The earliest appearance of these two non-regular triangulations 
that we know of is in  \cite{ConHen-moae}, although they
are closely related to Sch\"onhardt's classical 
example of a  non-convex
3-polytope that cannot be triangulated~\cite{Schoenhardt}.\footnote{We describe Sch\"onhardt's polyhedron and its relation to this example in Example~\ref{exm:moae-schoenhardt}.}

\end{example}

\begin{remark}
Suppose that two point sets 
$\A=\{a_1,\dots,a_n\}$ and $\B=\{b_1,\dots,b_n\}$ have
the same \emph{oriented matroid}~\cite{OM-book}, or \emph{order type}. 
This means that for every subset
$I\subset\{1,\dots,n\}$ of labels, the determinants of the point sets
$\{a_i:i\in I\}$ and $\{b_i:i\in I\}$ have the same sign.\footnote{Observe
that the  bijection between $\A$ and $\B$ implicit
by the labels is part of the defintion.} It is an easy exercise to check
that then $\A$ and $\B$ have \emph{the same} triangulations
and subdivisions.\footnote{More precisely, the implicit bijection between $\A$ and $\B$ induces a bijection between their
polyhedral subdivisions.}
However, they do not necessarily have the same \emph{regular} subdivisions.
For example, the points of example~\ref{exm:moae} are in general position and,
hence, their oriented matroid does not change by a small perturbation of coordinates. But any 
sufficiently generic perturbation makes one of the two non-regular triangulations $T_1$ and $T_2$ become regular.

Still,
the following is true~\cite{Santos-refine}: the \emph{existence} of non-regular triangulations
of $\A$  depends only on the oriented matroid of $\A$.
\end{remark}

\label{sub:secondary-gale}
\subsubsection*{The secondary polytope}

Let $L_\A$ denote the space of all lifting functions $w:\A\to \reals$
on a certain point set $\A\subset\reals^d$ with $n$ elements. In principle
$L_\A$ is isomorphic to $\reals^n$ in an obvious way; but we mod-out
functions that lift all of $\A$ to  a hyperplane, because adding
one of them to a given lifting function $w$ does not 
(combinatorially) change the lower envelope
of $\A_w$. We call these particular 
lifting functions \emph{affine}. They form a
linear subspace of dimension $d+1$ of $\reals^n$. Hence,
after we mod-out affine functions we have $L_\A\cong\reals^{n-d-1}$.

For a given polyhedral subdivision $T$ of
$\A$, the subset of $L_\A$ consisting of functions $w$ that produce
$T=T_w$, is a 
(relatively open) polyhedral cone;
that is, it is defined by a finite set of linear homogeneous equalities and strict inequalities. Equalities appear only if $T$ is not a triangulation
and express the fact that if $\sigma\in T$ is not affinely independent
 then $w$ must lift all $\sigma$ to lie in a hyperplane. Inequalities
express the fact that for each $\sigma\in T$ and point 
$a\in\A\setminus \sigma$, $a$ is lifted above the hyperplane spanned by
the lifting of $\sigma$. 

The polyhedral cones obtained for different choices 
of $T$ are glued together forming a polyhedral fan, that is, a ``cone over
a polyhedral complex'', called the \emph{secondary fan}
of $\A$.
The prototypical example of a fan is the normal
fan of a polytope, whose cones are the exterior normal cones of different 
faces of $P$. A seminal result in the theory of triangulations of
polytopes is that the secondary fan is actually polytopal; that is, it is the normal fan of a certain polytope:

\begin{theorem}[Gel'fand-Kapranov-Zelevinskii~\cite{GKZ-paper,GKZ-book}]
\label{thm:secondary}
For every point set $\A$ of $n$ points affinely spanning $\reals^d$
there is a polytope $\Sigma(\A)$ in $L_\A\cong\reals^{n-d-1}$
whose normal fan is the secondary fan of $\A$.
\end{theorem}

In particular, the poset of regular subdivisions of $\A$ is isomorphic
to the poset of faces of $\Sigma(\A)$. Vertices correspond to regular triangulations and $\Sigma(\A)$ itself (which is, by convention, 
considered a face) corresponds to the trivial subdivision.
The polytope $\Sigma(\A)$ is called the \emph{secondary polytope of $\A$}.

There are two standard ways to construct the secondary polytope
$\Sigma(\A)$ of a point set $\A$.%
\footnote{
Polytopality
of a fan is equivalent to 
the feasibility of a certain
system of linear equalities and  strict inequalities.
But here we mean more direct and intrinsic 
constructions of the secondary polytope.}
The original one, by Gel'fand, Kapranov and 
Zelevinskii~\cite{GKZ-paper,GKZ-book}
gives, for each regular triangulation $T$ of $\A$, coordinates
of the corresponding vertex $v_T$ of $\Sigma(\A)$ in terms of the 
volumes of simplices incident in $T$ to each point of $\A$.

The second one, by Billera and Sturmfels~\cite{BilStu},
 describes the whole polytope $\sigma(\A)$ as the
Minkowski integral of the fibers of the affine projection
 $\pi:\Delta_\A\to\conv(\A)$,
where $\Delta_\A$ is a simplex with $|\A|$ vertices (hence, of dimension
$|\A|-1$) and $\pi$  bijects the vertices of $\Delta_\A$ to $\A$~(see Theorem~{thm:fiber}).

\begin{example}(Example~\ref{exm:fivepoints} continued)
\label{exm:fivepoints-secondary}
Figure~\ref{fig:fivepoints-second} shows the secondary fan of
the five points. 
To mod-out affine functions we have taken
$w(a_1)=w(a_2)=w(a_3)=0$, and the horizontal and vertical coordinates
in the figure give the values of $w(a_4)$ and $w(a_5)$, respectively. 
The triangulation corresponding to each two-dimensional
cone is displayed.
\begin{figure}[htb]
\centerline{
\includegraphics[height=1.3 in]{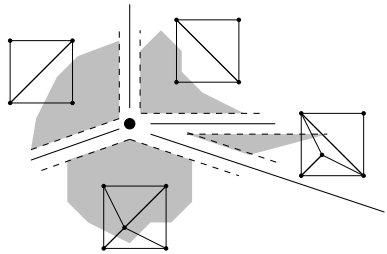}{}{}
}
\caption{The secondary fan of Example~\ref{exm:fivepoints}}
\label{fig:fivepoints-second}
\end{figure}
In this example all nine polyhedral subdivisions are regular
(in agreement with the  result of \cite{Lee-regular}
mentioned in Example~\ref{exm:moae})
and the secondary polytope is a quadrilateral.
\end{example}

\begin{example}(Example~\ref{exm:moae} continued)
\label{exm:moae-secondary}
The secondary polytope of this point set is 3-dimen\-sional, and contains
a hexagonal face corresponding to the regular subdivision
\[
\{1245,2356,1346,456\}.
\]
This regular subdivision can be refined to a triangulation in eight 
ways, by independently inserting a diagonal in the quadrilaterals
$1245$, $2356$ and $1346$. Six of these triangulations are regular,
and correspond to the vertices of the hexagon. The other two, $T_1$ and $T_2$,
are non-regular and they ``lie'' in the center of the hexagon.

We have mentioned that if the point set is perturbed slightly then one of the triangulations becomes regular. What happens in the secondary polytope 
is the following: the perturbation ``inflates'' the hexagon so that
the eight points (the
vertices of the hexagon and the two interior points representing $T_1$ and $T_2$)  become, combinatorially, the vertices of a cube. The points
corresponding to $T_1$ and $T_2$ move in opposite directions, one of them
going to the interior of the secondary polytope and the other becoming
a new vertex of it.
The hexagonal face gets refined into three 
quadrilaterals. Of course, the vertices of the hexagon also move in the
process, and are no longer coplanar.
\end{example}

\begin{example}(The convex $n$-gon and the associahedron)
All triangulations of a convex $n$-gon are regular and their
number is the $n-2$nd Catalan
number
\[
C_{n-2}:=\frac{1}{n-1}\left(\genfrac{}{}{0pt}{}{2n-4}{ n-2}\right).
\]
The corresponding secondary polytope is called the \emph{associahedron}.
The name comes from the fact that there is a bijection between 
triangulations of the $n$-gon and the ways to put the $n-2$ parentheses
in an associative product of $n-1$ factors. 

The associahedron
is a classical object in combinatorics,
first studied\footnote{As a combinatorial cell complex, without an explicit
polytopal realization.}
by Stasheff and Tamari \cite{Tamari,Stasheff}. It was shown to 
be polytopal by Haiman (unpublished) and Lee~\cite{Lee}.
That its diameter equals $2n-10$ ``for 
every sufficiently big $n$''\footnote{%
Sleator~et~al.~do not say ``how big'' is ``sufficiently big''
in their statement, 
but conjecture that $n\ge 13$ is enough. We consider
this an interesting and somehow shameful open question.}
was shown by Sleator, Tarjan and Thurston~\cite{SlTaTh}, with
motivations coming from theoretical computer science and
tools from  hyperbolic geometry.
\end{example}

\begin{remark}
\label{rem:vectorconfs}
In Sections~\ref{sub:topological-combinatorics} 
and~\ref{sub:algebraic-geometry} we will mention
triangulations of a set of \emph{vectors} rather than \emph{points}. They are defined
exactly as triangulations of point sets, just changing the word
``affinely'' to ``linearly'' and the operator ``$\conv$'' to ``$\pos$'' (``positive span'')
in Definition~\ref{defi:triangulation}. Put differently,
a triangulation of a vector set $\A\subset\reals^{d+1}$ 
is a simplicial fan covering
$\pos(\A)$ and whose rays are in the positive directions of (not
necessarily all) the elements of $\A$. Equivalently, and perhaps
closer to readers familiar with classical geometry, we can, without
loss of generality, normalize all vectors of $\A$ to lie in the unit sphere $S^d$.
Then, triangulations of $\A$ are
the geodesic triangulations, with vertices contained in $\A$, of the spherical
convex hull of $\A$.

The existence and properties of regular subdivisions and secondary
fans (and of the bistellar flips introduced in the next section)
generalize almost without change to vector configurations.\footnote{
Although with one notable difference. 
For a general vector configuration
not every function $w:\A\to \reals$ produces a lift with a well-defined
``lower envelope''. Only the functions that do, namely those for
which a linear hyperplane exists containing or lying 
below all the lifted vectors,
define a regular polyhedral subdivision. These functions
form a cone in $L_\A$. The secondary fan is still 
well-defined but, of course, it cannot
be the 
%
normal fan of a polytope.
It is, however, the normal fan of an unbounded
convex polyhedron, called the \emph{secondary polyhedron} of $\A$~\cite{BiGeSt}.
}
\end{remark}

\subsection{Geometric bistellar flips}
\label{sub:bistellar}

\subsubsection*{Flips as polyhedral subdivisions}
In order to introduce the notion of local move (flip) between triangulations
of $\A$, we use the secondary fan as a guiding light: whatever
our definition is, restricted to regular triangulations a flip
should correspond to crossing
%
%
a ``wall'' between two full-dimensional cones in the  
secondary fan; that is, a flip 
between two regular triangulations $T_1$ and $T_2$ can be regarded as
certain regular subdivision 
$T_0$ with the property that its
only two regular refinements are precisely $T_1$ and $T_2$.
Some thought
will convince the reader that the necessary and sufficient condition for a lifting function $w:\A\to\reals$ to produce a $T_w$
with this property is that \emph{there is a unique 
minimal affinely dependent
subset in $\A$ whose lifting is contained in some lower facet
of the lifted point set $\A_w$}. 
%
This leads to the following simple, although perhaps not very practical, definition.

\begin{definition}
\label{defi:flip-subdivision}
Let $T$ be a (not-necessarily regular) subdivision of a point set $\A$.
We say that $T$ \emph{is a flip} 
if there is a unique affinely dependent 
subset $C\in \A$ contained in some cell of $T$.
\end{definition}


\begin{lemma}
\label{lemma:flip-subdivision}
If $T$ is a flip, then there are exactly two proper refinements of $T$,
which are both triangulations.
\end{lemma}

\begin{proof}
Let $T_1$ be a refinement of $T$. Let $C$ be the unique affinely dependent
subset of $\A$ contained in some cell of $T$. Each cell of $T$ containing $C$ gets refined in $T_1$, while each cell
 not
containing $C$ is also a cell in $T_1$. 

The statement
then follows from the understanding of the combinatorics of point sets
with a unique affinely dependent subset $C$. 
Let $S$ be such a point set. Each point
in $S\setminus C$ is affinely independent of the rest, so $S$ is
an ``iterated cone'' over $C$. In particular, there is a face $F$
of $S$ such that $S\cap F=C$ and every refinement of $S$
consists of a refinement of $F$ coned to the points of $S\setminus C$.
Moreover, all cells of $T$ containing $C$ must have $F$ refined the same way,
so that there is a bijection between the refinements of $T$ 
and the polyhedral subdivisions of $C$, as a point set. The result then
follows from the fact (see below) that a minimal affinely dependent
set $C$ has exactly three subdivisions: the trivial one and two triangulations.
\end{proof}

This lemma allows us to understand a flip, even in the non-regular case, 
as a relation or a transformation between its two refinements. 
This is the usual usage of the word ``flip'', and our next topic.

\subsubsection*{Flips as elementary changes}
A minimal affinely dependent set $C$ is called a \emph{circuit}
in geometric combinatorics. The points in a circuit $C=\{c_1,\dots,c_k\}$
satisfy a unique (up to a constant) affine dependence equation
$\lambda_1c_1+\cdots+\lambda_kc_k=0$ with $\sum \lambda_i=0$, and
all the $\lambda_i$ must be non zero (or otherwise $C$ is not 
minimally dependent).
This affine dependence implicitly decomposes $C$ into two subsets
\[
C_+=\{c_i : \lambda_i>0\},\qquad
C_-=\{c_i : \lambda_i<0\}.
\]
The pair $(C_+,C_-)$ is usually called a \emph{signed} or
\emph{oriented} circuit. We will slightly abuse notation 
and speak of ``the circuit $C=(C_+,C_-)$'', unless we need
to emphasize the distinction between the set $C$ 
(the \emph{support} of the circuit)
and its partition.

A more geometric description
 is that $(C_+,C_-)$ is the only partition of $C$ into
two subsets whose convex hulls intersect, and that 
they intersect in their relative interiors. This is usually called
\emph{Radon's property}~\cite{Radon} and the oriented circuit a
\emph{Radon partition}.

Spanning and affinely independent subsets
of $C$ are all the sets of the form $C\setminus \{c_i\}$. Moreover,
by Radon's property two such sets $C\setminus \{c_i\}$ and $C\setminus \{c_j\}$ can be cells in the same triangulation of $C$ 
if and only if $c_i$ and $c_j$
lie in the same side of the partition. In other words:

\begin{lemma}
A circuit $C=(C_+,C_-)$ has exactly two
triangulations:
\[
T^C_+:=\{C\setminus\{c_i\} : c_i\in C_+\}, \qquad
T^C_-:=\{ C\setminus\{c_i\} : c_i\in C_-\}.
\]
\end{lemma}

This leads to a second definition of flip, equivalent to 
Definition~\ref{defi:flip-subdivision}, but more 
operational. This is the definition originally devised
by Gel'fand, Kapranov and Zelevinskii~\cite{GKZ-book}.
The \emph{link}
of a set $\tau\subseteq\A$ in a triangulation $T$ of $\A$
is defined as 
\[
\link_T(\tau):= \{\rho\subseteq \A: \rho\cap\tau=\emptyset,
\rho\cup\tau\in T\}.
\]

\begin{definition}
\label{defi:flip-operation}
Let $T_1$ be a triangulation of a point set $\A$. Suppose
that $T_1$ contains one of the triangulations, say $T^C_+$,
of a circuit $C=(C_+,C_-)$. Suppose also
that all cells $\tau\in T^C_+$ have the same link
in $T_1$, and call it $L$.

Then, we say that $C$ \emph{supports a geometric bistellar flip}
(or a \emph{flip}, for short)  in $T_1$ and that
the following triangulation $T_2$ of $\A$ is obtained from $T_1$ by
this flip:
\[
T_2:= T_1 
\setminus \{\rho\cap\tau: \rho\in L, \tau\in T^C_+\}
\cup \{\rho\cap\tau: \rho\in L, \tau\in T^C_-\}.
\]

If $i=|C_+|$ and $j=|C_-|$ we say that the flip is of type $(i,j)$.
Flips of types $(1,j)$ and $(i,1)$ are called,
\emph{insertion} and \emph{deletion} flips, since they add or
remove a vertex in the triangulation.

The \emph{graph of flips} of $\A$ has as vertices all the 
triangulations of $\A$ and as edges the geometric bistellar flips between them.
\end{definition}

Of course, an $(i,j)$ flip can always be reversed, giving a $(j,i)$ flip.
The reason for the words ``geometric bistellar'' in our flips
can be found in Section~\ref{sub:combinatorial-topology}.

\begin{example}(Examples~\ref{exm:fivepoints} and~\ref{exm:fivepoints-secondary}
continued)
\label{exm:fivepoints-flips}
The change between the two top triangulations in Figure~\ref{fig:fivepoints-second} is a $(2,2)$ flip, as is the change between the two bottom ones.
The flip from the top-right to the bottom-right is a $(1,3)$ flip
(``1 triangle disappears and 3 are inserted'') and the flip from
the top-left to the bottom-left is a $(1,2)$ flip (``one edge is 
removed, together with its link,  and two are inserted, with the same link'').
The latter is
supported in the circuit formed by the three collinear points.
\end{example}

We omit the proof of the  following natural statement.
\begin{theorem}
\label{thm:flips-equivalence}
Definitions~\ref{defi:flip-subdivision} and~\ref{defi:flip-operation}
are equivalent: two triangulations $T_1$ and $T_2$ of a point set $\A$
are connceted by a flip in the sense of~\ref{defi:flip-operation}
if and only if they are the two proper refinements of a flip in
the sense of~\ref{defi:flip-subdivision}.
\end{theorem}

The following two facts are proved in~\cite{Santos-refine}:
\begin{remark}
\label{rem:good-reasons}
\begin{enumerate}
\item  
If all proper refinements of a subdivision $T$ are 
triangulations, then $T$ has
exactly two of them and $T$ is a flip. That is to say,
flips are exactly the ``next-to-minimal'' elements in the refinement
poset of all subdivisions of $\A$.

\item Every non-regular subdivision can be refined
to a non-regular triangulation~\cite{Santos-noflips}. In particular, 
not only edges of the secondary polytope correspond to flips between two
regular triangulations, but also every flip between two regular 
triangulations corresponds to an edge.
\end{enumerate}
\end{remark}

\subsubsection*{Detecting flips}
Definitions~\ref{defi:flip-subdivision} and~\ref{defi:flip-operation}
are both based on the existence of a \emph{flippable circuit} $C$ 
with certain properties. But in order to detect flips only some circuits
need to be checked:

\begin{lemma}
\label{lemma:walls}
Every flip in a triangulation $T$ other than an insertion flip is supported
in a circuit contained in the union of two adjacent cells
of $T$.
\end{lemma}

Observe that the circuit contained in two adjacent cells always exists
and is unique.
Also, that the insertion flips left aside in this statement
are easy
to detect:\footnote{We mean, theoretically. 
Algorithmically, insertion flips are far from trivial
since they imply locating the simplex of $T$ that contains the
point $a$ to be inserted, which takes about
the logarithm of the number of simplices in $T$. This is 
very expensive, since algorithms
in computational geometry that use flipping in triangulations
usually are designed to take constant time per flip other than
an insertion flip. See Section~\ref{sub:computational-geometry}.}
There is one for each point $a\in\A$ not used in $T$, that inserts
the point $a$ by subdividing the minimum (perhaps not full-dimensional)
simplex $\tau\subseteq\sigma\in T$ such 
that $a\in\conv(\tau)$. The
flippable circuit is $(\{a\},\tau)$.

\begin{proof}
Let $C=(C_+,C_-)$ be a circuit that supports a flip in $T$,
with $|C_+|\ge 2$. Observe that $|C_+|$ is also the number of 
many maximal simplices in $T^C_+$, so let 
$tau_1$ and $\tau_2$ be two of them, which differ
in a single element, and let $\rho$ be an element of $\link_T(\tau_1)=
\link_T(\tau_2)$. Then, $\rho\cup\tau_1$ and  $\rho\cup\tau_2$
are adjacent cells in $T$ and $C$ is the unique
circuit contained in $\tau_1\cup\tau_2\cup\rho$.
\end{proof}

\subsubsection*{Monotone sequences of flips}
The graph of flips among regular triangulations of 
a point set $\A$ of dimension $d$ is connected,
since it is the graph of a polytope.\footnote{
Even more, it is $(|\A|-d-1)$-connected.
Remember that a graph is called $k$-connected if removing
less than $k$ vertices from it it stays connected. A classical theorem
of Balinski~\cite{Ziegler-book} says that the graph of a $k$-polytope
is $k$-connected.}
A fundamental fact exploited in computational geometry is that
one can actually flip between regular triangulations
\emph{monotonically}, in the following sense.

Let $w:\A\to \reals$ be a certain generic lifting
function. We can use $w$ to lift
every triangulation $T$ of $\A$ as a function 
$H_{T,w}:\conv(\A)\to\reals$,
by affinely interpolating $w$ in each cell of $T$. We say that
$T_1 <_w T_2$ (``$T_1$ is below $T_2$, with respect to $w$'') if  $H_{T_1,w}\le H_{T_2,w}$ pointwise and $H_{T_1,w}\ne H_{T_2,w}$ globally. This defines a partial order $<_w$ 
on the set of all triangulations, whose global minimum and maximum are, respectively, $T_w$ and $T_{-w}$.\footnote{In case they are triangulations. If not, every triangulation that refines $T_w$ or $T_{-w}$ is, respectively, minimal or maximal.}

\begin{definition}
\label{defi:monotone}
A sequence of flips is monotone with respect to $w$ if every flip
goes from a triangulation $T$ to a triangulation $T' <_w T$.
\end{definition}

By definition of the secondary polytope $\Sigma(\A)$
as having the secondary fan as its normal fan, 
lifting functions are linear functionals on it. 
Then, it is no surprise that for the regular triangulations $T_1$ and
$T_2$ corresponding to vertices $v_{T_1}$ and $v_{T_2}$ of
the secondary polytope one has\footnote{The same is true for non-regular
triangulations. The point $v_T$ is well-defined, via
the Gel'fand-Kapranov-Zelevinskii coordinates for the secondary polytope,
even if $T$ is not regular. The only difference is that if $T$ is not regular then
 $v_T$ is not a vertex of the secondary polytope.}:
\[
T_1 <_w T_2 \qquad \Rightarrow \qquad
\langle w, v_{T_1} \rangle < \langle w, v_{T_1} \rangle.
\]
In fact, $\langle w, v_{T} \rangle$ equals the volume between
the graphs of the functions $H_{T_w,w}$ and $H_{T,w}$.
Since the converse implication holds whenever $T_1$ and $T_2$
are related by a flip, we have:

\begin{lemma}
\label{lemma:monotone}
For every lifting function $w$ and every regular triangulation $T$
there is a $w$-monotone sequence of flips from $T$ to the regular
triangulation $T_w$.
\end{lemma}

If $T$ is not regular this may be false, even in 
dimension 2:
\begin{example} (Examples~\ref{exm:moae}
 and~\ref{exm:moae-secondary} continued)
\label{exm:moae-schoenhardt}
Let $\A$ be the point configuration of Example~\ref{exm:moae}
(see Figure~\ref{fig:moae-triangs}),
except perturbed by 
slightly rotating the interior triangle 
``$123$'' counter-clockwise. That is,
{\small
\[
\A=
\bordermatrix{
& {a_1} & {a_2} & {a_3} & {a_4} & {a_5} & {a_6} \cr
&4-\varepsilon & 0 & \varepsilon & 2 & 1 & 1 \cr
&\varepsilon & 4-\varepsilon & 0 & 1 & 2 & 1 \cr
&0 & \varepsilon & 4-\varepsilon & 1 & 1 & 2 \cr
},
\]
}
for a small $\varepsilon >0$.
This perturbation
keeps the triangulation $T_1$ non-regular and makes $T_2$ regular.
Let $w:\A\to\reals$ lift the exterior triangle $123$
to height zero and the interior
triangle $456$ to height one. The graph of $H_{T_2,w}$ is 
a strictly
concave surface (that is, $T_2=T_{-w}$) and
there is no $w$-monotone flip in $T_1$, since its only three 
flips are the diagonal-edge flips on
``$16$'', ``$24$'' and ``$35$'', which are ``towards $H_{T_2,w}$''.
This example appeared in~\cite{EdeSha}.

Another explanation of why no $w$-monotone flip exists in $T_1$ is that 
when we close the graph of the function
$H_{T_1,w}$ by adding to it the triangle $123$, 
it becomes  a non-convex polyhedron $P$ with the property that no tetrahedron 
(with vertices contained in those of $P$)
is completely contained in the region enclosed by $P$. 
This polyhedron is affinely equivalent to 
Sch\"onhardt's~\cite{Schoenhardt} classical
example of a non-convex polyhedron in $\reals^3$
that cannot be triangulated without
additional vertices.%
\end{example}

\subsection{The cases of small dimension or few points}
\label{sub:state-of-the-art}


Throughout this section $\A$ denotes a point set with $n$ elements
and dimension $d$.

\paragraph{Sets with few points.}
If $n=d+1$,
then $\A$ is independent and the trivial subdivision is its unique triangulation. If $n=d+2$ then $\A$ has a unique circuit and
 exactly two triangulations, connected by a flip. 
If  $n=d+3$, it was proved by Lee~\cite{Lee-regular}
that all triangulations are regular. Since the secondary fan 
is 2-dimensional,
the secondary polytope is a polygon, whose graph (a cycle) is 
the graph of flips.
If $n=d+4$, then $\A$ can
have non-regular triangulations (see Example~\ref{exm:moae}). Still,
it is proven in~\cite{AzaolaSantos} that every triangulation has 
at least three flips and that the flip-graph is 3-connected.

For point sets with $n=d+5$ the flip-graph is not known to be always 
connected.

\paragraph{Dimension 1.}
Triangulating a one-dimensional point set
is just choosing which of the interior points
are  used. That is, $n$ points in dimension 1 have $2^{n-2}$
triangulations. The flip-graph is the graph of an $(n-2)$-dimensional
cube and all triangulations are regular. 
The secondary polytope is the same cube.

\paragraph{Dimension 2.}
In dimension two the graph of 
\emph{$(2,2)$-flips among triangulations
using all points of $\A$}~\footnote{This is the graph usually considered
in two-dimensional computational geometry literature.} is known to
be connected since long~\cite{Lawson}, and connectivity of the
whole graph---including the triangulations that do not use all points and the insertion or deletion flips---is
straightforward from that. 
%
Even more, one can always flip 
monotonically\footnote{With respect to the lift $w(a):=||a||^2$.}
from any triangulation to the Delaunay triangulation using only
$(2,2)$ flips. Quadratically many (with respect to the number of points) flips are sometimes necessary and always suffice
(see, e.g., ~\cite[p. 11]{Edels-book}).

However, with general flips:

\begin{proposition}
The flip-graph of
any $\A\subseteq\reals^2$ has diameter smaller than $4n$.
\end{proposition}

\begin{proof}
Let $a$ be an extremal point of $\A$ and
$T$ an arbitrary triangulation. If $T$ has triangles 
not incident to $a$ then there is at least
a flip that decreases the number of them (proof left to the reader). 
Since the number of triangles in a planar triangulation with $v_i$
interior verticess and $v_b$ boundary vertices is exactly $2v_i+v_b-2$ (by Euler's 
formula) we can flip from any triangulation to one with every triangle
incident to $a$ in at most $2v_i+v_b-3<2n - n_b$ flips. 

Now, eaxctly as in 
the 1-dimensional case, the graph of flips between
triangulations in which every triangle is incident to $a$ 
is the graph of a cube of dimension equal to
the number of ``boundary but non-extremal'' points of $\A$.
%
Hence, we can flip between
any two triangulations in $(2n - n_b) + n_b + (2n - n_b)< 4n$ flips.
\end{proof}

\begin{remark}
\label{rem:lexicographic}
The preceding proof is another example of monotone flipping,
this time with respect to any
lifting function $w:\A\to \reals$ with $w(a)<<w(b)$, 
for all $b\in\A\setminus \{a\}$. In essence, this lifting produces
the so-called \emph{pulling} triangulation of $\A$. More precisely,
for a point set $\A$ in arbitrary dimension and a given 
ordering $a_1,\dots,a_n$ of the points in $\A$ one 
defines~\cite{OM-book,Lee-regular,Lee-survey,Ziegler-book}:
\begin{itemize}
\item The \emph{pulling triangulation} of $\A$, as the regular triangulation
given by the lift $w(a_i):=-t^i$, for a sufficiently big constant $t\in\reals$.
It can be recursively constructed as the
triangulation that joins the last point $a_n$ to the pulling
triangulation of every facet of $\conv(\A)$ that does not contain $a_n$.

\item The \emph{pushing triangulation} of $\A$, as the regular triangulation
given by the lift $w(a_i):=t^i$, for a sufficiently big constant $t\in\reals$.
It can be recursively constructed as the
triangulation that contains  $T_{n-1}$ and joins $a_n$ to the part of the boundary of  $T_{n-1}$ visible from $a_n$,
where $T_{n-1}$  is the pushing triangulation
 of $\A\setminus\{a_n\}$. 
\end{itemize}
Pushing and pulling triangulations are examples of \emph{lexicographic}
triangula\-tions, defined by the lifts $w(a_i):=\pm t^i$ for sufficiently
big $t$. 
\end{remark}

Summing up,  monotone flipping in the plane:
(a) Works even for
non-regular triangulations if the ``objective function'' $w$
is either the Delaunay or a lexicographic one
(the proof for the pushing case is left to the reader).
%
(b) It gives a linear sequence of flips for the pulling case,
but may produce a quadratic one for the Delaunay case.
%
(c) It does not work for arbitrary $w$ 
(Example~\ref{exm:moae-schoenhardt}).

Let us also mention that
in dimension two every triangulation
is known to have at least $n-3$ flips~\cite{LoSaUr}
(the dimension of the secondary polytope), and at least $\lceil (n-4)/2 \rceil$ 
of them of type $(2,2)$~\cite{HuNoUr}.
The flip-graph is not known to be $(n-3)$-connected.

\paragraph{Dimension 3.}
Things start to get complicated:

If $\A$ is in convex position~\footnote{\emph{Convex position} means that all points are vertices of $\conv(\A)$.} then 
every triangulation of it has at least $n-4$
flips~\cite{LoSaUr}, but otherwise $\A$ can have
 flip-deficient triangulations.\footnote{We say a triangulation
is  \emph{flip-deficient} if it has less than $n-d-1$ flips; that
is, less than the dimension of the secondary polytope.}
The smallest possible example, with eight points, 
is described in~\cite{AzaolaSantos}, based on
 Example~\ref{exm:moae}.
%
Actually, for every $n$
there are triangulations with essentially $n^2$ vertices and only 
$O(n)$ flips~\cite{Santos-fewflips}. This is true even 
in general position.\footnote{Although this
is not mentioned in~\cite{Santos-fewflips}, the construction there can be
perturbed without a significant addition of flips.}
The flip-graph  is not known to be connected,
even if $\A$ is in convex and general position.
The main obstacle to proving  connectivity
(in case it holds!) is probably that
one cannot, in general, \emph{monotonically}
flip to either the Delaunay,
the pushing, or the pulling triangulations. For the Delaunay triangulation
this was shown in~\cite{Joe-stucked}. For the other two we describe here an example. 
%

\begin{example}(Examples~\ref{exm:moae},~\ref{exm:moae-secondary} and~\ref{exm:moae-schoenhardt} continued)
\label{exm:moae-schoenhardt-3d}
Let $\A$ consist of the following eight points in dimension three: 
{\small
\[
\A=
\bordermatrix{
& {a_1} & {a_2} & {a_3} & {b_1} & {b_2} & {b_3} & {c_1} & {c_2} \cr
&4-\varepsilon & 0 & \varepsilon & 2 & 1 & 1 & {4}/{3} & {4}/{3}  \cr
&\varepsilon & 4-\varepsilon & 0 & 1 & 2 & 1 & {4}/{3} & {4}/{3}  \cr
&0 & \varepsilon & 4-\varepsilon & 1 & 1 & 2 & {4}/{3} & {4}/{3}  \cr
& 0 & 0 & 0 & 1 & 1 & 1 & 10 & -10 \cr
},
\]
}
The first six points are exactly the (lifted) point set of 
Example~\ref{exm:moae-schoenhardt}, and have the property that
no tetrahedron with vertices contained in these six points is
contained in the non-convex Sch\"onhardt polyhedron $P$
having as boundary triangles $\{a_1,a_2,a_3\}$,
$\{b_1,b_2,b_3\}$, $\{a_i,a_{i+1},b_{i}\}$ and $\{a_{i+1},b_{i},b_{i+1}\}$
(the latter for the three values of $i$, and with indices regarded modulo three).
The last two points $c_1$ and $c_2$ of the configuration lie far above
and far below this polyhedron. $c_1$ sees every face of $P$ except the 
big triangle $\{a_1,a_2,a_3\}$, while $c_2$ sees only this triangle.

Let $T$ be  the triangulation $T$ of $\A$ obtained 
removing the big triangle from the boundary of $P$, 
and joining the other seven triangles to both $c_1$ and $c_2$.
We leave it to the reader to check that there is no monotone sequence
of flips towards the pushing triangulation with respect to any ordering
ending in $c_2$%
, and there is no monotone sequence
of flips towards the pulling triangulation with respect to any ordering
ending in $c_1$.
%
%

%
%
%
%
\end{example}

\paragraph{Higher dimension.}
There are the following known examples of ``bad behavior'':

\begin{itemize}
\item In \textbf{dimension four}, there are triangulations with arbitrarily many vertices and a bounded number of flips~\cite{Santos-fewflips}.
They are constructed adding several layers of ``the same'' triangulated
3-sphere one after another.

\item In \textbf{dimension five}, there are point sets with a disconnected
graph of triangulations~\cite{Santos-hilbert}. The smallest one
known has 26 points
, but one with 50 points is easier to describe:
It is the Cartesian  product of $\{0,1\}$ with the vertex set and the
centroid of a regular 24-cell.

\item In \textbf{dimension six}, there are triangulations 
without flips at all~\cite{Santos-noflips}. The example is again a certesian product, now of a very
simple configuration of four points in $\reals^2$ and a not-so-simple
(although related to the 24-cell too) configuration of 81 points 
in $\reals^4$.
There are also point sets \emph{in general position} with a disconnected
graph of triangulations 
(Section~\ref{sec:construction} of this paper and~\cite{Santos-17}). Only 17 points are needed.
\end{itemize}


\section{The context}
\label{sec:context}
\subsection{Bistellar flips and computational geometry}
\label{sub:computational-geometry}
The first and most frequently considered flips in computational
geometry are $(2,2)$ flips in 2-dimensional point sets. A
 seminal result of Lawson~\cite{Lawson} says that every
triangulation can be monotonically transformed in the Delaunay
triangulation by a sequence of $O(n^2)$ such flips. 

Around 
1990,\footnote{Birth year of secondary
polytopes and geometric bistellar flips as we have defined them~\cite{GKZ-paper}.}
B.~Joe generalizes flips to three dimensions, defining the 
$(3,2)$ and $(2,3)$ flips (and, implicitly, also the ``insertion'' $(1,4)$
flips). He realizes that one cannot, in general, monotonically flip from 
any triangulation to the Delaunay triangulation~\cite{Joe-stucked} but,
still, the following incremental algorithm works: 
insert the points one by one, each 
by an insertion flip in the Delaunay triangulation of the already
inserted points.
After each insertion, monotonically
flip to the new Delaunay triangulation by
flips that increase the star of the inserted point.

V.~T.~Rajan~\cite{Rajan} does essentially the same in arbitrary dimension  and
 Edelsbrunner and Shah~\cite{EdeSha}, already
aware at least partially of the theory of secondary polytopes,
generalize this to flipping towards the regular triangulation $T_w$
by $w$-monotone flips, for an arbitrary $w$.

If one disregards the efficiency of the algorithm, the main
result of \cite{EdeSha} follows easily from Lemma~\ref{lemma:monotone}.
But efficiency is the main point in computational geometry,
and one of the important features in \cite{Joe-good} and \cite{EdeSha} is to 
show that the sequence of flips can be found and performed
spending constant time per flip (in fixed dimension).
%
An exception to this time bound
are the insertion steps.~
Theoretically, they are 
just another case of flip. But in the algorithm they have
a totally different role since they involve locating where the
new point needs to be inserted. 
To get good time bounds for the location step,
the standard incremental algorithm is ``randomized'',\footnote{
That is, the ordering in which the points are inserted is
considered random among the $n!$ possible orderings.
This trick was first introduced by Guibas, Knuth and Sharir~\cite{GuiKnuSha}
in Joe's two-dimensional incremental algorithm.} and
it is proved that  the total 
\emph{expected} time taken by the $n$ insertion steps is 
bounded above  by $O(n\log n)$ 
in the plane and $O(n^{\lceil d/2 \rceil})$ in higher dimension.
The latter is the same as the worst-case size of the Delaunay
triangulation, or actually of any triangulation.

This incremental-randomized-flipping method can be considered the standard 
algorithm for the  Delaunay triangulation
in current computational geometry. {It is the only one described in
the textbooks~\cite{deBerg-etal} and \cite{Edels-book}.
In the survey~\cite{AurenKlein}, it is 
the first of four described in the plane
but the only one detailed in dimension three, as ``the most intuitive
and easy to implement''.}

\begin{remark}
\label{rem:generalposition}
Computational geometry literature normally 
only considers full-dimen\-sional flips;
that is, flips of type $(i,j)$ with $i+j=d+2$.
In particular,~\cite{AurenKlein},~\cite{deBerg-etal},~\cite{Edels-book},~\cite{EdeSha} and~\cite{Joe-good}
describe the incremental flipping algorithm only 
for point sets in general position.
%
The only mention in those references
 to the effect of allowing special position
in the flipping process seems to be that, according to~\cite{Edels-book},
for the six vertices of a regular octahedron
\emph{``none of the three tetrahedrizations permits the application
of a two-to-three or a three-to-two flip. The flip graph thus
consists of three isolated nodes''}.
%

However, with the general definition of flip
the incremental-flipping algorithm can be directly applied
to point sets in special position, as done recently
by Shewchunk~\cite{Shewchunk}.
Shewchunk's algorithm actually  computes the 
so-called  \emph{constrained regular triangulation} of the point set
for any lift $w$ and constrain complex $K$.
This is defined as the unique\footnote{If it exists, which is not always the case.}
triangulation $T$ 
containing $K$ and in which every simplex of $T\setminus K$
is lifted by $w$ to have a locally convex star.%
\footnote{Shewchunk's algorithm is incremental, treating 
the simplices in $K$
similarly to the points in the standard incremental algorithm:
they are inserted one by one (in increasing
order of dimension) and after each insertion the
regular, constrained to the already added simplices,
triangulation is updated using 
geometric bistellar flips.  The algorithm's  running time
is $O(n^{\lfloor d/2 \rfloor+1}\log n)$.
The extra $\log n$ factor comes from a 
priority queue that is needed to decide in which order the
flips are performed, to make sure that no ``local optima''
instead of the true constrained Delaunay triangulation,
is reached. The extra $n$ factor (only in even dimension) 
is what randomization saves in the standard incremental-flipping
algorithm, but it is unclear to us whether 
randomization would do the same here.
}

%
\end{remark}

%

\subsection{Bistellar flips and combinatorial topology}
\label{sub:combinatorial-topology}

Bistellar flips can be defined at a purely combinatorial level, for an abstract simplicial complex.
Let $\Delta$ be a simplicial complex, and let $\sigma\in\Delta$ be a simplex, of any dimension.
The \emph{stellar subdivision} on the simplex $\sigma$ is the simplicial complex obtained inserting
a point in the relative interior of $\sigma$. This subdivides $\sigma$, and every simplex $\tau$ containing it, into
$\dim{\sigma} +1$ simplices of the same dimension. Two simplicial complexes $\Delta_1$ and $\Delta_2$
are said to differ in a \emph{bistellar flip} if there are simplices $\sigma_1\in \Delta_1$ and $\sigma_2\in \Delta_2$
such that the stellar subdivisions of $\Delta_1$ and $\Delta_2$ on them produce the same simplicial complex.
The bistellar operation from $\Delta_1$ to $\Delta_2$ is said to be of type $(i,j)$ if
$i=\dim{\sigma_1} +1$ and $j=\dim{\sigma_2} +1$.
%
Observe that geometric bistellar flips, as defined in Definition~\ref{defi:flip-operation}, are
combinatorially bistellar flips.

Combinatorial bistellar flips have been proposed as an algorithmic tool for exploring the space of
triangulations of a manifold
\footnote{Besides its intrinsic interest,
this problem arises in quantum gravity modelization~\cite{ACM,Nabutovsky}.}
or to recognize the topological type of a simplicial manifold~\cite{BjornerLutz,Lickorish}.
 In particular, Pachner~\cite{Pachner} has shown that any
two triangulations of PL-homeomorphic manifolds are connected by a sequence of topological bistellar flips.
But for this connectivity result
 additional vertices are allowed to be inserted into the complex, via
flips of type $(i,1)$. 

The situation is much different if we do not allow insertion flips: 
Dougherty et al.~\cite{DoFaMu} show that there is a 
topological triangulation of the 3-sphere, with 15 vertices,
that does not admit any flip other than insertion flips.%
\footnote{Dougherty et al. only say that their 
triangulation does not have any
$(3,2)$, $(2,3)$ or $(1,4)$ flips, which are the ``full-dimensional'' types of flips.
 But their arguments
prove that even considering degenerate flips, the only possible ones in
their triangulation 
are insertion flips of type $(i,1)$. Indeed, the two basic properties that
their triangulations has are that (a) its graph is complete, which prevents
flips of type $(3,2)$, but also $(2,2)$ and $(1,2)$ and
(b) no edge is incident to exactly three tetrahedra, which prevents
flips of type $(1,4)$ and $(2,3)$, but also $(1,3)$.
}
If this triangulation was realizable
gometrically in $\reals^3$ (removing from the 3-sphere
the interior of any particular simplex) it would provide a triangulation in dimension three without
any geometric bistellar flips.
Unfortunately, Dougherty et al. show that it cannot be geometrically embedded.


\subsection{Bistellar flips and topological combinatorics}
\label{sub:topological-combinatorics}
A standard construction in topological combinatorics~\cite{Bjorner}
is to associate to a poset $P$ its \emph{order complex}:
an abstract simplicial complex
whose vertices are the elements of $P$ and whose simplices are the
finite chains (totally ordered subsets) of $P$. In this sense
one can speak of the topology of the poset. If the poset has
a unique maximum (as is the case 
with the refinement poset of subdivisions of a point set $\A$)
or minimum,
one usually removes them or otherwise the order complex is trivially contractible
(that is, homotopy equivalent to a point).
This is what we mean when we say that the
refinement poset of subdivisions of the 
point set of Section~\ref{sec:construction} is not connected.

The refinement poset of polyhedral subdivisions
of $\A$ is usually called the \emph{Baues poset} of $\A$ and its study
is \emph{the generalized Baues problem}.
To be precise, Baues posets were introduced 
implicitly in \cite{BilStu} and explicitly in \cite{BiKaSt} in a 
more general situation where one has an affine 
projection $\pi$ from the vertex set of a polytope $P\in\reals^{d'}$
to a lower dimensional affine space $\reals^d$. In this general setting,
one considers the point set $\A:=\pi(\vertices(P))$
and is interested in the polyhedral subdivisions 
of $\A$ that are compatible
with $\pi$ in a certain sense (basically, that the preimage of every
cell is the set of vertices of a face in $P$).
In the special case were $P$ is a simplex (and hence $d'=n-1$, where
$n$ is the number of points in $\A$) every polyhedral subdivision is 
compatible. This is the case of primal interest in this paper,
but there are at least the following two other cases that have attracted attention.
(See \cite{Reiner-survey} for a
very complete account of different contexts in which Baues posets appear,
and~\cite[Chapter 9]{Ziegler-book} for a different treatment of the topic):

\begin{itemize}
\item 
When $P$ is a cube, its projection is a {\em
zonotope} $Z$ and the $\pi$-compatible subdivisions are the {\em zonotopal tilings} of $Z$~\cite{Ziegler-book}. The finest ones are \emph{cubical tilings}, related by \emph{cubical flips}.

\item When $d=1$ and $P$ is arbitrary, the $\pi$-compatible subdivisions
are called \emph{cellular strings}, since they correspond to monotone
sequences of faces of $P$. The finest ones are \emph{monotone paths} of edges and are related by \emph{polygon flips}.
\end{itemize}

The name \emph{Baues} for these posets comes from the fact that H.~J.~Baues
was interested in their homotopy type in a very particular
case (in which, among other things, $d=1$) and conjectured it
to be that of a sphere of dimension $d'-2$~\cite{Baues}.
Billera et al.~\cite{BiKaSt} proved this conjecture for all Baues posets with $d=1$, and the conjecture that the same happened for arbitrary $d$ 
(with the dimension of the sphere being now $d'-d-1$)
became known as the~\emph{generalized Baues conjecture}.
It was inspired by the fact that the \emph{fiber polytope}
associated to the projection $\pi$---a generalization of the secondary
polytope, introduced in~\cite{BilStu}---has dimension $d'-d$ and its face lattice is
naturally embedded in the Baues poset. 

Even after the conjecture in its full generality was disproved by 
a relatively simple example with $d'=5$ and $d=2$~\cite{RamZie},
the cases where $P$ is either a simplex or a cube remained of interest.
As we have said, the latter is disproved in the present paper
for the first time. The former remains open and has connections to 
oriented matroid theory, as we now show.


\medskip

Recall that the oriented matroid (or order type)
of a point set $\A$ 
of dimension $d$ (or of a vector configuration
of rank $d+1$) is just the 
information contained in the map $\genfrac{(}{)}{0pt}{}{\A}{ d+1}\to\{-1,0,+1\}$
that associates to each $(d+1)$-element subset of $\A$ the sign of 
its determinant (that is, its orientation).
But oriented matroids (see~\cite{OM-book} as a general reference)
are axiomatically defined structures
which may or may not  be realizable as the oriented matroids
of a real configuration, in much the same way as, for example, 
a topological space may or may not be metrizable. 

It turns out that the theory
of triangulations of point and vector configurations generalizes
nicely to the context of perhaps-non-realizable oriented
matroids, with the role of 
regular triangulations by
the so-called \emph{lifting triangulations}: triangulations
that can be defined by an oriented matroid lift
(see \cite[Section 9.6]{OM-book}
or \cite[Section 4]{Santos-OMtri}).

One of the basic facts in oriented matroid duality is that the lifts
of an oriented matroid $\M$ are in bijection to the one-point
extensions of its dual $\M^*$. 
In particular, the
space of lifts of $\M$ 
equals the so-called extension space
of the dual oriented matroid $\M^*$. 
Here, both the space of lifts and the space of extensions are defined
as the simplicial complexes associated to the natural poset structures in the
set of all lifts/extensions of the oriented matroid.
%
%
This makes the following
conjecture of Sturmfels and Ziegler~\cite{StuZie} be relevant to
this paper:
\begin{conjecture}
The extension space of a realizable oriented matroid 
of rank $r$ is homotopy equivalent to a sphere of dimension $r-1$.
\end{conjecture}

The  reader may be surprised that we call this a 
conjecture: if the extension space of an oriented matroid is the
analogue of a secondary fan, shouldn't the
extension space of a realizable oriented matroid be automatically
``a fan'', hence a sphere? Well, no: even if an oriented matroid
$\M$ is realizable, some of its extensions may not be realizable.
Those will appear in the extension space.
Even worse, if $\M$ is realized as a vector configuration
$\A$, some realizable extensions of $\M$ may only be realizable as
extensions of other realizations of $\M$.
Actually, Sturmfels and Ziegler show that the space of realizable extensions
of a realizable oriented matroid \emph{does not} in general have the homotopy 
type of a sphere!

\begin{example} (Example~\ref{exm:moae} continued)
\label{exm:moae-om}
Consider the point configuration of Example~\ref{exm:moae} (two parallel
triangles one inside the other). An additional point added
to this configuration represents an extension of the
underlying oriented matroid. In particular, there is an extension
by a point collinear with each of the three pairs of corresponding
vertices of the two triangles.

But any small perturbation of the point set gives another realization
of the same oriented matroid, since the original point set is in general position. However,
this perturbation will, in general, not keep the lines through
those three pairs of vertices colliding. So, the extension we have described
is no longer realizable as a geometric extension of the new realization.

\end{example}

%
There is a class of configurations
specially interesting in this context:
the so-called Lawrence polytopes. A \emph{Lawrence oriented
matroid} is an oriented matroid whose dual is centrally symmetric.
Similarly, a \emph{Lawrence polytope} is a polytope whose vertex set
has a centrally symmetric \emph{Gale transform}. There is essentially one
Lawrence polytope associated to each and every realizable oriented matroid.
The following result
is a combination of a theorem of Bohne and Dress (see~\cite{Ziegler-book}, 
for example) and one of the author of this paper~\cite{HuRaSa,Santos-OMtri}:

\begin{theorem}
\label{thm:lawrence}
Let $\M$ be a realizable oriented matroid and let 
$P$ be the associated Lawrence polytope. Then, the following three posets
are isomorphic:
\begin{enumerate}
\item The refinement poset of polyhedral subdivisions of $P$.
\item The extension space of the (also realizable) dual oriented matroid
$\M^*$.
\item The refinement poset of zonotopal tilings of the zonotope
associated to (any realization of) $\M$.
 \end{enumerate}
\end{theorem}

\begin{corollary}
\label{coro:lawrence}
The following three statements are equivalent:
\begin{enumerate}
\item The generalized Baues conjecture for the polyhedral subdivisions of 
Lawrence polytopes.
\item The extension space conjecture for realizable oriented matroids.
\item The generalized Baues conjecture for the zonotopal tilings of 
zonotopes.
\end{enumerate}
\end{corollary}

Moreover, if $\A$ is a point configuration and $P$ its associated 
Lawrence polytope, then there is a surjective map between
the poset of subdivisions of $P$ and the poset of \emph{lifting}
(in the oriented matroid sense) subdivisions of $\A$.
This follows from that facts that ``Lawrence polytopes only have lifting subdivisions'' and ``lifting subdivisions can be lifted to the Lawrence polytope'', both proved in~\cite{Santos-OMtri}. 

In particular, if the flip-graph of a certain point set $\A$
is not connected and has
lifting triangulations in several connected components,
then the graph of cubical flips between zonotopal tilings
of a certain zonotope is not connected either,
thus answering question 1.3 in \cite{Reiner-survey}.
If, moreover, $\A$ is in general
position, it would disprove the three statements in
Corollary~\ref{coro:lawrence}.
We do not know whether the disconnected flip-graph
in Section~\ref{sec:construction} has this property.
The examples in~\cite{Santos-noflips,Santos-hilbert}
are easily seen to be based in non-lifting triangulations.

\begin{remark}
The extension space conjecture is the case
$k=d-1$ of the following far-reaching conjecture by MacPherson, Mn\"ev and
Ziegler \cite[Conjecture 11]{Reiner-survey}: that the poset of all
strong images of rank $k$ of any  realizable oriented matroid $\M$ of rank
$d$ (the so-called {\em OM-Grassmannian of rank $k$} of $\M$)
is homotopy equivalent to the real Grassmannian 
$G^k(\reals^d)$. This conjecture is relevant in
matroid bundle theory \cite{Anderson} and the 
{\em combinatorial differential geometry} of MacPherson
\cite{MacPherson}.

 An important achievement
in this context is the recent result of Biss~\cite{Biss}
proving this conjecture whenever $\M$ is a ``free oriented matroid''.
In this case the OM-Grassmannian is  the space of all
oriented matroids of a given cardinality and rank, usually called
the MacPhersonian.
The result of Biss includes the case $n=\infty$ (in which the MacPhersonian
is defined as a direct limit of all the MacPhersonians of a given rank)
and implies that the theory of ``oriented matroid bundles for
combinatorial differential manifolds'' developed by MacPherson~\cite{MacPherson} is equivalent to the theory of real
vector bundles on real differential manifolds.
A first, seminal, result in this direction was the
``combinatorial formula'' by Gel'fand and MacPherson for the 
Pontrjagin class of a triangulated manifold~\cite{GelfandMac}.
\end{remark}

\subsection{Bistellar flips and algebraic geometry}
\label{sub:algebraic-geometry}

Bistellar flips are related to algebraic geometry from their very birth. Indeed, 
Definition~\ref{defi:flip-operation},
as well as that of secondary polytope and Theroem~\ref{thm:secondary} were first given
by Gel'fand, Kapranov and Zelevinskii during their study of discriminants of a sparse polynomial~\cite{GKZ-paper}.
By a sparse polynomial we mean, here, a multivariate polynomial $f$ whose coefficients are 
considered parameters but whose set of (exponent vectors of) monomials is a fixed 
 point set $\A\subseteq\integers^{d}$. Gel'fand, Kapranov and Zelevinskii prove that the secondary 
 polytope of $\A$ equals the Newton polytope of the Chow polynomial of $f$, where the Chow polynomial 
 is a certain resultant defined in terms of $f$. Similarly, the secondary polytope is related
 to the discriminant of $f$ (the $\A$-discriminant) although a bit less directly: it is a Minkowski
 summand of the Newton polytope of the $\A$-discriminant.

A stronger, and more classical, relation between triangulations of point sets and algebraic geometry comes
from the theory of toric varieties~\cite{Fulton,Oda}.
As is well-known, every rational convex polyhedral fan $\Sigma$ (in our language, every 
polyhedral
subdivision of a rational vector configuration) has an associated toric variety $X_\Sigma$, of the same dimension. 
$X_\Sigma$ is non-singular if and only if $\Sigma$ is simplicial (i.e., a triangulation) and unimodular. The
latter means that every cone is spanned by integer vectors with determinant $\pm 1$.
If $\Sigma$ is a non-unimodular triangulation, then $X_\Sigma$ is an orbifold; that is, it has only quotient singularities.

A stellar subdivision, that is, an insertion flip, in $\Sigma$ corresponds to an equivariant blow-up in 
$X_\Sigma$. Hence, a deletion flip produces a blow-down and a general flip produces a 
blow-up followed by a blow down. In this sense, the connectivity question for triangulations
of a vector configuration is closely related to the following result, conjectured by Oda~\cite{MiyakeOda} and proved by Morelli~\cite{Morelli}.\footnote{Morelli's original paper contained some minor errors. See~\cite{AbMaRa} and~\cite{Morelli2} for corrections. The weak form of Oda's conjecture was independently proved by
Wlodarczyk~\cite{Wlodarczyk}.}

\begin{theorem}
Every proper and equivariant birational map $f : X_\Sigma\to X_{\Sigma'}$  between two nonsingular toric varieties can be factorized into a sequence of blowups and blowdowns with centers being smooth closed orbits
(weak Oda's conjecture). Moreover, we can insist on the sequence to consist of first a sequence of only blowups and then one of only blowdowns (strong Oda's conjecture). 
\end{theorem}

More precisely, Oda's conjecture, in its weak form, is equivalent to saying that every pair of unimodular simplicial fans
can be connected by a sequence of bistellar flips passing only through unimodular fans (and, actually, it is proved this way). But observe that
in this result the set of vectors allowed to be used is not fixed in advance: additional ones are allowed to be flipped-in and eventually flipped-out. Our construction in~\cite{Santos-hilbert} actually shows that the result is not true if we don't allow for extra vectors to be inserted.

\medskip

The relation of the graph of flips to toric geometry is even closer if one looks at certain schemes associated to 
a toric variety. In order to define them 
we first look at secondary polytopes in a different way, as
a particular case of fiber polytopes~\cite{BilStu}: 

Assume that $\A$ is an integer point configuration and 
let $\Delta$ 
be the unit simplex of dimension $|\A|-1$ in $\reals^{|\A|}$. Let $Q=\conv(\A)$ and let 
$\pi:\Delta\to Q$ be the affine projection sending the vertices of $\Delta$ to
$\A$. 
The  chamber complex of $\A$ is the coarsest common refinement of
all its triangulations. It is a polyhedral complex with the
property that for any $b$ and $b'$ in the same {\em chamber} the
fibers $\pi^{-1}(b)$ and $\pi^{-1}(b')$ are polytopes with
the same normal fan. 

\begin{theorem}[Billera et al.~\cite{BilStu}]
\label{thm:fiber}
The secondary polytope of $\A$ equals the
Min\-kowski integral of $\pi^{-1}(b)$ over $Q$.
\end{theorem}

Combinatorially, then, the secondary polytope of $\A$ 
equals the Minkowski sum of a finite number of $\pi^{-1}(b)$'s,
with one $b$ chosen in each chamber.

Now, for each $b\in Q$, consider the toric variety associated to the 
normal fan of the fiber $\pi^{-1}(b)$. Since the normal fan is the same whenever 
$b$ and $b'$ lie in (the relative interior of) the same cell of the chamber complex, 
we denote this toric variety $V_\sigma$, where $\sigma$ is a cell 
(of any dimension) of the chamber complex.
If $b\in\sigma$ and  $b'\in\tau$ for two
chambers with $\tau\subseteq \overline{\sigma}$ then the normal fan of
$\pi^{-1}(b)$ refines the normal fan of $\pi^{-1}(b')$, 
which implies that there is a
natural equivariant morphism $f_{\sigma\tau}: V_{\sigma}\to V_{\tau}$. 
We finally denote $\Lambda_\A:=\invlim V_\sigma$ 
the inverse limit of all the $V_\sigma$ and
morphisms $V_{\sigma\tau}$. It has the following two interpretations:

\begin{enumerate}
\item Let $X_\Delta$ be the projective space of dimension $|\A|-1$, 
which is the toric variety associated with the simplex $\Delta$ (what
follows is valid for any polytope $\Delta$).
The toric varieties $V_{\sigma}$ are
the different toric geometric invariant theory quotients of 
$X_\Delta$ modulo the algebraic sub-torus whose
characters are the monomials with exponents in $\A$
\cite[Section~3]{KaStZe}.
$\Lambda_\A$ is the inverse limit of all of
them, which contains the Chow quotient as an irreducible component
\cite[Section~4]{KaStZe}.

\item In \cite{Alexeev}, Alexeev is interested, among other things,  in 
the moduli space $M$ of {\em stable semi-abelic toric pairs}
for an integer polytope $Q$
(see Sections 1.1.A and 1.2.B in
\cite{Alexeev} for the definitions).
The author shows that there is a
finite morphism $M\to \Lambda_\A$ (Corollary 2.11.11),
where $\A$ is the set of all integer points in $Q$, 
and uses
$\Lambda_\A$ (that he denotes $M_{simp}$)
as a simplified model for studying $M$.

Although there $\A$ is assumed
to be the set of all lattice points in a polytope, the connection of
$\Lambda_\A$ with $\Sigma_c(\A)$
carried out in the proof of the following theorem is independent of
this fact. 
\end{enumerate}

\begin{theorem}
\label{thm:alexeev}
The scheme $\Lambda_\A$ is connected if and only if the graph of triangulations
of $\A$ is connected.
\end{theorem}

\begin{proof}[Proof (Sketch)\ :]
Alexeev introduces the following poset structure on the set of all polyhedral subdivisions of $\A$:
Given two subdivisions $S_1$ and $S_2$ we consider $S_1<S_2$ if: 
(a) $S_1$ refines $S_2$,
(b) the restriction of $S_1$ to each cell $B$ of $S_2$ is a regular subdivision $S_B$ of $B$, and
(c) the lifting functions of the regular subdivisions of cells of $S_2$ can be chosen
so that the restrictions of them to common faces of cells differ by an affine function.

This poset is called the ``coherent poset of subdivisions of $\A$'' in \cite{Santos-noflips}, to distinguish it from
the usual poset of subdivisions, where only the first condition (refinement) is imposed.
Then, he shows that the scheme $\Lambda_\A$ is connected if and only if the coherent refinement poset
is connected. (More precisely, he shows that there is a natural moment map defined on $\Lambda_\A$
whose image is the topological model of the poset). In turn, it is proven in~\cite{Santos-noflips} that 
the coherent refinement poset is connected if and only if the graph of triangulations of $\A$ is connected.
\end{proof}


A second scheme that relates  triangulations and toric geometry
is precisely the so-called toric Hilbert scheme. The toric ideal 
$I_\A\subseteq K[x_1,\dots,x_n]$ associated to $\A=\{a_1,\dots,a_n\}\in\reals^d$ is generated by the binomials 
\[
\{x^\lambda-x^\mu : \lambda,\mu \in \naturals^n, \sum \lambda_i a_i = \sum \mu_i a_i\}.
\]
Here, $x^\lambda:=x_1^{\lambda_1}\cdots x_n^{\lambda_n}$. In other words, $I_\A$ is the lattice ideal
of the lattice of integer affine dependences among $\A$.
$\A$ defines the following $\A$-grading of monomials in $K[x_1,\dots,x_n]$: the $A$-degree of
$x^\lambda$ is the vector $x_1^{\lambda_1}\cdots x_n^{\lambda_n}\in \integers^{d}$. Of course, $I_\A$ 
is homogeneous with respect to this grading.

If $I$ is another $\A$-homogeneous ideal, the Hilbert function of $I$ is the map $\integers^d\to \naturals$
defined by $b\mapsto \dim_{K} I_b$ where $I_b$ is the part of $I$ of degree $b$. The \emph{toric Hilbert scheme} of $\A$ consists, as a set, of all the $\A$-homogenous ideals with the same Hilbert function as the
toric ideal $I_\A$. It contains $I_\A$ as well as all its initial ideals, which form an irreducible component in its 
scheme structure.

The toric Hilbert scheme was introduced by Sturmfels in~\cite{Sturmfels-Agraded} (see also \cite{Sturmfels-book}) although its scheme structure was explicited later by Peeava and Stillman \cite{PeeSti},
who ask whether non-connected toric Hilbert schemes exist.

Sturmfels shows, among other things, that there is a natural map from the toric Hilbert scheme to 
the set of polyhedral subdivisions of $\A$. Moreover, the map is continuous when the latter is
given either the poset topology or the ``coherent poset topology'' introduced in the proof of Theorem~\ref{thm:alexeev}. The map is not surjective in general, so disconnected graphs of triangulations do not
automatically imply disconnected Hilbert schemes.\footnote{Haiman and Sturmfels~\cite{HaiStu}
have shown that this map factors as a morphism  from the toric Hilbert scheme to the
scheme $\Lambda_\A$ of the previous discussion, followed by the natural map from that scheme
to the poset of subdivisions. The first map is the non-surjective one.}
However, 
Maclagan and Thomas~\cite{MacTho}, modifying the arguments of Theorem~\ref{thm:alexeev}, show
that the image of the map contains at least al the unimodular triangulations of $\A$. In particular:

\begin{corollary}
If the graph of triangulations of an integer point configuration $\A$ is not connected and contains unimodular
triangulations in non-regular connected components, then the toric Hilbert scheme of $\A$ is not connected.
\end{corollary}

The example in~\cite{Santos-hilbert} satisfies the hypothesis of this corollary. Hence:

\begin{theorem}[Santos~\cite{Santos-hilbert}]
Let $\A_{50}\subset\reals^5$ be the point set $\A_{25}\times\{0,1\}$ where $\A_{25}\subset\reals^4$ consists of the centroid and the 24 vertices of a regular 24-cell.
The toric Hilbert scheme of $\A$ and the scheme $\Lambda_A$ defined above are both non-connected.
They have at least 13 connected components, each with at least  $3^{48}$ torus-fixed points.
\end{theorem}

\section{A construction}
\label{sec:construction}

Let $\A(t)\subset\reals^6$ be the point set defined by the columns of the following 
matrix, where $t$ denotes a positive real number. 
The matrix is written in two pieces for typographic reasons.
As usual, the first row is just a homogenization coordinate:
{\small
\begin{align*}
\A(t):=
\bordermatrix
{
& O & a_1^+(t) & a_2^+(t) & a_3^+(t) & a_4^+(t) & a_5^+(t) & a_6^+(t) & a_7^+(t) & a_8^+(t) \cr
& 1  & 1 & 1 & 1 & 1 & 1 & 1 & 1 & 1\cr 
x_1& 0  & 1 &-t & 0 & 0 & 1 & t & 0 & 0 \cr
x_2& 0  & t & 1 & 0 & 0 &-t & 1 & 0 & 0 \cr
x_3& 0  & 0 & 0 & 1 &-t & 0 & 0 & 1 & t \cr
x_4& 0  & 0 & 0 & t & 1 & 0 & 0 &-t & 1 \cr
x_5& 0  & \sqrt{2} & 1 & 0 & -1 & -\sqrt{2} & -1 & 0 & 1 \cr 
x_6& 0  & 0 & 1 & \sqrt{2} & 1 & 0 & -1 & -\sqrt{2} & -1  \cr
}
\\
\dots
\bordermatrix
{
& a_1^-(t) & a_2^-(t) & a_3^-(t) & a_4^-(t) & a_5^-(t) & a_6^-(t) & a_7^-(t) & a_8^-(t) \cr
  & 1 & 1 & 1 & 1 & 1 & 1 & 1 & 1\cr 
  &-1 & t & 0 & 0 &-1 &-t & 0 & 0 \cr
  &-t &-1 & 0 & 0 & t &-1 & 0 & 0 \cr
  & 0 & 0 &-1 & t & 0 & 0 &-1 &-t \cr
  & 0 & 0 &-t &-1 & 0 & 0 & t &-1 \cr
  & \sqrt{2} & 1 & 0 & -1 & -\sqrt{2} & -1 & 0 & 1 \cr 
  & 0 & 1 & \sqrt{2} & 1 & 0 & -1 & -\sqrt{2} & -1  \cr
}.
\end{align*}
}
$\A(t)$ is
not in general position.
For example, for every $i=1,2,3,4$ we have:
\[
a_i^+(t)+a_{i+4}^+(t)+a_i^-(t)+a_{i+4}^-(t) = 4 O.
\]
However, it is ``sufficiently in general position'' for the following
to be true:

\begin{theorem}
\label{thm:17}
If $t$ is sufficiently small and $\A'(t)$ is any perturbation of $\A(t)$ in general position, then the graph of triangulations of $\A'(t)$ is not connected.
\end{theorem}

When we say that a point set $\A'$ is a perturbation of another one $\A$ 
with the same cardinal $n$ and dimension $d$ we mean that
all the determinants of $d+1$ points that are not zero in $\A$ keep their
sign in $\A'$.\footnote{In oriented matroid language, the oriented matroid of
$\A$ is a weak image of that of $\A'$.} This concept also allows us to be
precise as to how small do we need $t$ to be. Any $t$ such that $\A(t)$
is a perturbation of $\A(0)$ works. 

\medskip

The proof of Theorem~\ref{thm:17} will appear in~\cite{Santos-17}.
Here we only give a description of the combinatorics of $\A(t)$ and
the ingredients that make the proof work.
We look at $\A(0)$ first. In it:
\begin{itemize}
\item The projection to the first four coordinates $x_1,\dots,x_4$ sends 
the eight pairs of points $\{a_i^+(t),a_{i+4}^+(t)\}$, and $\{a_i^-(t),a_{i+4}^-(t)\}$ ($i=1,2,3,4$) to the eight vertices of a 4-dimensional cross-polytope (that is,
to the standard basis vectors and their opposites).
\item The projection to the last two coordinates $x_5,x_6$ sends 
the eight pairs of points $\{a_i^+(t),a_{i}^-(t)\}$ ($i=1,\dots,8$) to the eight 
vertices of a regular octagon.
\end{itemize}

The configuration $\A(0)$ already has a disconnected graph of triangulations. 

\begin{theorem}
\label{thm:17-degenerate}
There is a triangulation $K$ of the boundary of $\conv(\A(0))$ with the following two properties:
\begin{enumerate}
\item There are triangulations of $\A(0)$ inducing $K$ on the boundary.
\item No flip in a triangulation of $\A(0)$ inducing $K$ on the boundary 
affects the boundary.
\end{enumerate}
\end{theorem}

In fact, there are eight such triangulations. Hence:

\begin{corollary}
The flip-graph of $\A(0)$ has at least nine connected components.%
\footnote{
Here, the ninth component is the one containing all the regular triangulations.}
\end{corollary}

Of course, to describe the triangulation $K$ of the boundary of $\conv(\A_0)$ we
need only specify how we triangulate each  non-simplicial facet.
The facets of $\conv(\A(0))$ are 96 simplices, and 16 non-simplicial 
facets $F_{\delta_1,\delta_2,\delta_3,\delta_4}$ ($\delta_i\in\{+,-\}$), each with eight vertices.
More precisely,
\[
F_{\delta_1,\delta_2,\delta_3,\delta_4}=
\{
a_1^{\delta_1}(0), a_2^{\delta_2}(0), a_3^{\delta_3}(0), a_4^{\delta_4}(0),
a_5^{\delta_1}(0), a_6^{\delta_2}(0), a_7^{\delta_3}(0), a_8^{\delta_4}(0)
\}.
\]
All the $F_{*,*,*,*}$'s are equivalent under affine symmetries of $\A(0)$.
For example, they are transitively permuted by the sixteen sign 
changes on the first four coordinates. 
Hence, the crucial point
in the proof of Theorem~\ref{thm:17-degenerate} is to understand the
triangulations of the point set $F_{+,+,+,+}$. 
This point set has
 dimension $d=5$  and only eight ($=d+3$) points. In particular, all its
triangulations are regular and their graph of flips is a cycle.
Moreover, it is easy to check\footnote{For example, noting that a Gale transform of
$F_{+,+,+,+}$ consists again of the eight vertices of a regular octagon, except in different order.} that:

\begin{lemma}
\label{lemma:17}
\begin{enumerate}
\item $\conv(F_{+,+,+,+})$ has 12 facets. 
Eight of them are simplices and the other four have six points
each, forming a $(3,3)$ circuit.
In particular,
there are sixteen ways to triangulate the boundary of $F_{+,+,+,+}$.
\item $F_{+,+,+,+}$ has eight triangulations. 
\item Each 
flip in a triangulation of $F_{+,+,+,+}$ keeps the triangulation
induced in three of the non-simplicial facets and switches the 
triangulation in the other.
\end{enumerate}
\end{lemma}




To construct the complex $K$ of Theorem~\ref{thm:17-degenerate} we choose the triangulations of the individual
$F_{*,*,*,*}$ such that for every non-simplicial facet $G$ of an $F_{\delta_1,\delta_2,\delta_3,\delta_4}$, the triangulations 
chosen on $F_{\delta_1,\delta_2,\delta_3,\delta_4}$ and on the neighbor
$F_{\delta'_1,\delta'_2,\delta'_3,\delta'_4}$ 
agree on $G$ and one of them has the property that
no flip on it changes the triangulation induced in $G$.
In these conditions, no flip in any of the triangulations of the
$F_{*,*,*,*}$'s is possible, since it would be incompatible with the triangulation of one of its neighbors.

\begin{example}
\label{exm:prism}
Lemma~\ref{lemma:17} implies, in particular, that only eight
of the sixteen triangulations of the boundary of 
$F_{+,+,+,+}$ can be extended to the interior
(without using additional
interior points as vertices).
Similar behavior occurs also
in three-dimensional examples such as the 
set of vertices of a cube or a triangular prism.

Let us analyze the latter. 
It has three non-simplicial facets, whose vertex sets are
$(2,2)$ circuits;
in particular,
there are eight ways to triangulate its boundary.
But only six of them extend to the interior (all except the
two ``cyclic'' ones). Each 
flip in a triangulation of $F_{+,+,+,+}$ keeps the triangulation
induced in two of the non-simplicial facets and switches the 
triangulation in the other one.%
\footnote{
The reader probably has noticed the similarities between this
example and the configuration $F_{+,+,+,+}$.
These similarities, and the fact that the constructions
in \cite{Santos-noflips} and \cite{Santos-hilbert} are 
ultimately based on glueing
triangular prisms to one another, reflect the truth in 
(an instance of)
Gian Carlo Rota's fifth lesson~\cite{Rota}.
}
\end{example}

Let us now look at the perturbations $\A(t)$ and $\A'(t)$.
The fact that $\A(t)$ (or $\A'(t)$) is a perturbation of $\A(0)$ implies
that every triangulation of $\A(0)$ is  still a geometric simplicial complex 
on $\A(t)$, except it may not cover the whole
convex hull. In particular, the triangulation $K$ of the boundary
of $\conv(\A_(0))$ mentioned in Theorem~\ref{thm:17-degenerate} can be
embedded as a simplicial complex on $\A(t)$. We still call $K$ this
perturbed simplicial complex.
Then, Theorem~\ref{thm:17}
follows from the following more precise statement.

\begin{theorem}
\label{thm:17-perturbed}
Let $t$ be a sufficiently small and positive constant. Then,
\begin{enumerate}
\item There are triangulations of $\A(t)$ containing the simplicial complex $K$.
\item If $T$ is a triangulation of $\A(t)$ containing the simplicial complex $K$, then every triangulation obtained from $T$ by a flip contains the simplicial complex $K$. In particular, the graph of triangulations of 
$\A(t)$ is not connected.
\item The previous two statements remain true if $\A(t)$ is perturbed into
general position in an arbitrary way.
\end{enumerate}
\end{theorem}

\end{document}